\documentclass[a4paper,11pt]{amsart}

\usepackage[a4paper, left=2.3cm, right=2.3cm, top=2cm, bottom=2cm]{geometry}
\usepackage{meintex}
\usepackage{caption}
\usepackage{subcaption}
\usepackage[algo2e,norelsize]{algorithm2e}
\usepackage{algpseudocode}
\usepackage{mdframed}
\usepackage[section]{placeins}

\title[On the Numerics of Mean-field Optimal Control]{Instantaneous control of interacting particle systems in the  mean-field limit }

\author[M.~Burger]{Martin Burger}
\address[M.~Burger]{Friedrich-Alexander Universität Erlangen-Nürnberg}
\email{martin.burger@fau.de}

\author[R.~Pinnau]{Ren\'e Pinnau}
\address[R.~Pinnau]{TU Kaiserslautern}
\email{pinnau@mathematik.uni-kl.de}

\author[C.~Totzeck]{Claudia Totzeck}
\address[C.~Totzeck]{TU Kaiserslautern}
\email[Corresponding author]{totzeck@mathematik.uni-kl.de}

\author[O.~Tse]{Oliver Tse}
\address[O.~Tse]{Technische Universiteit Eindhoven}
\email{o.t.c.tse@tue.nl}

\author[A.~Roth]{Andreas Roth}
\address[A.~Roth]{TU Kaiserslautern}
\email{roth@mathematik.uni-kl.de}

\begin{document}
\maketitle
\begin{abstract}
Controlling large particle systems in collective dynamics by a few agents is a subject of high practical importance, e.g., in evacuation dynamics.  In this paper we study an instantaneous  control approach to steer an  interacting particle system into a certain spatial region by repulsive forces from a few external agents, which might be interpreted as shepherd dogs leading sheep to their home. We introduce an appropriate mathematical model and the corresponding optimization problem. In particular, we are interested in the interaction of numerous particles, which can be approximated by a mean-field equation. Due to the high-dimensional phase space this will require a tailored optimization strategy. The arising control problems are solved using adjoint information to compute the descent directions.  Numerical results on the microscopic and the macroscopic level indicate the convergence of optimal controls and optimal states in the mean-field limit,i.e., for an increasing number of particles.
\end{abstract}
\section{Introduction}
In the last decades, the behavior of large particle systems and their mean-field limits were intensively investigated on both the theoretical and computational level \cite{Golse,AlbiPareschi,CarrilloMills,carrillo2009double}. 
Large groups of individuals like flocks of birds and schools of fish, and their attractive and repulsive interaction were considered, which lead to models of different types of collective behavior such as flocking or milling, and a thorough study of their stability, see, e.g., \cite{reynolds1987flocks,czirok2000collective,barbaro2009discrete,maini2001mathematical} or \cite{CarrilloMills} for a detailed overview.
The models were refined using vision cones, self-propulsion and orientation alignment of neighbors\cite{Dorsogna,CarrilloMills,yates2011refining}. As the behavior of large particle groups of the same type was well understood, the interest in self-organized systems interacting with few external agents arose. This concept and its kinetic limit was first numerically investigated by Albi and Pareschi in \cite{AlbiPareschi}. Using interaction potentials introduced by Cucker-Smale \cite{cucker2007emergent} or D'Orsogna et al \cite{Dorsogna}, they showed numerically, that the collective behavior of large groups coincides with the behavior of the kinetic model as the number of individuals tends to infinity. While the models for the interaction of particles became more and more realistic, the complexity was significantly increased. At first, interaction potentials were chosen smooth in order to have well-defined derivatives \cite{cucker2007emergent,cucker2007mathematics,Dorsogna}. Later, it was even shown that special classes of singular interaction potentials allow passing to the mean-field limit  \cite{hauray2007n,bolley2011stochastic}. Nowadays, particle games are also employed in the field of global optimization \cite{CBO,CBO2}.

The Model Predictive Control approach, in particular, instantaneous control was applied successfully to traffic flow problems \cite{HertyKlar} and the Navier-Stokes Equations \cite{HinzeNavier}. In the present work, we use an instantaneous approach to control the interaction of a huge crowd of individuals with a few external agents. The focus is here on the appropriate mathematical models and the tailored numerical approaches to solve these large scale optimization problems. In particular, we are interested in the behavior of microscopic and macroscopic optimization problems in the  corresponding mean-field  limit.The control parameters are the velocities of the external agents. The cost functionals are designed in such a way that the external agents aim to lead the crowd to a predefined destination while penalizing their kinetic energy. Exemplarily, one might consider a crowd of sheep guided by dogs. 

A similar sparse optimization problem was analytically  investigated  in \cite{FornasierSolombrino}, where the existence of optimal controls on the microscopic level is shown and the concept of $\Gamma$-convergence is used to perform the mean-field limit in a sparse optimal control setting.

Numerically, the computational effort of the state solutions alone is very challenging due to the high-dimensional phase space. On the one hand, we have the pairwise interactions in the microscopic system and on the other hand the mean-field equations are of Vlasov type, which yields a four dimensional problem in two spatial dimensions. The complexity of the optimization problems is even worse, since the state systems need to be solved several times and a huge amount of data needs to be stored. Hence, there is no hope of solving these problems using black-box optimization approaches. Instead, tailored optimization algorithms are required, which are typically based on derivative information \cite{Pinnau}. 

Recently, many researchers have investigated related problems: in \cite{invisible} an evacuation problem using invisible agents in a crowd was studied and a similar scenario with visible agents in~\cite{Evacuation}, \cite{Borzi} considers a Fokker-Planck feedback control strategy for crowd motion and \cite{zuazua}, where optimal strategies for driving a mobile agent in a guidance by repulsion model are studied. Optimal control problems for transport processes are discussed, e.g.,  in  \cite{buttazzo}. An extensive overview on the actual mathematical approaches for behavioral social systems can be found in \cite{review_bellomo}.


Here, we are going to construct appropriate numerical approaches  which allow us to verify computationally  that the controls of the microscopic problem converge to the optimal control of the kinetic problem as the number of individuals increases. The numerical algorithms are  based on first-order derivative information using the adjoint variables for  interacting particle system and  the mean-field equation, respectively.

The manuscript is organized as follows. First, we describe the microscopic interacting particle system and the corresponding kinetic mean-field equation in Section 2. The cost functionals for the  constrained optimization problems are introduced in Section 3, where we also specify the instantaneous control approach we use for the numerical simulations.  In Section 4 the associated first-order optimality conditions are derived for both the microscopic and the macroscopic optimization problems. The numerical schemes for the state and adjoint equations are presented in Section 5, where also the instantaneous control algorithm is introduced. Numerical results underlining our approach and confirming the desired behavior in the mean-field limit are discussed in Section 6. Concluding remarks are given in Section 7.
\section{Microscopic and Mean-Field Optimal Control Problems}\label{sec:state}
First, we describe the agent-based model and its corresponding mean-field limit, which are going to be the respective state systems for the control problems considered later on.

\subsection{Microscopic Model}
Let $D\ge 1$ denote the dimension of the spatial and the velocity space. The considered particle system consists of $N \in \mathbb{N}$ particles of the same type and $M$ external agents with $N\gg M$. Let $ I=[0,T]$ be the time interval of observation. Then, the particles and agents are represented by state vectors
\begin{equation*}
 x_i,v_i,d_m,u_m \colon  I \rightarrow \mathbb{R}^D, \qquad \text{for } i = 1,\dots N \text{ and }m=1,\dots,M.
\end{equation*}
The vectors 
\begin{equation*}
 \x(t) = (x_i(t))_{i=1,\dots,N}\quad \text{ and } \quad \ve(t) = (v_i(t))_{i=1,\dots,N}, 
\end{equation*}
 denote the positions and velocities of the particles and 
 \begin{equation*}
 \dog(t) = (d_m(t))_{m=1,\dots,M}\quad \text{ and } \quad \con(t) = (u_m(t))_{m=1,\dots,M},
 \end{equation*}
 the positions and velocities of the external agents, respectively (compare  also \cite{AlbiPareschi}). 
 
 \begin{rem}
 
 Since the time dependence is clear, we often write $\x$ or 
$x_i$ instead of $\x(t)$ and $x_i(t)$,  respectively.  Note, that we write the vectors $\x$ and $\ve$ in bold to make clear when denoting the positions and the velocities of the individuals and when we refer to the position and velocity space with variables $x$ and $v$ in the mean-field setting.
\end{rem}

The interactions of the individuals and external agents are modeled using potentials $\Phi_j$, $j=1,2$ which satisfy the following assumption
\begin{itemize} 
 \item [\textbf{(A)}]  The potentials $\Phi_j\colon\mathbb{R}^D\to\mathbb{R}$ are radially symmetric and continuously differentiable, with $\nabla \Phi_j$ locally Lipschitz and globally bounded.
\end{itemize}
\begin{rem}
These rather strict assumptions allow for existence and uniqueness of solutions to the adjoint and state systems. The latter are essential for defining the reduced cost functional needed in the algorithms for numerical investigation.
\end{rem}
To simplify the presentation, we denote for any $x,y\in\mathbb{R}^D$ the interaction forces 
\[
 K_j(x,y) = (\nabla\Phi_j)(x-y),\qquad j=1,2.
\]
Since $\Phi_j$ is symmetric, we have that $K_j(x,y)=-K_j(y,x)$ for all $x,y\in\mathbb{R}^D$.

Altogether, this leads to the particle system
\begin{subequations}\label{longODEsys}
 \begin{align}
 \frac{\dd}{\dd t}x_i &= v_i,\qquad \frac{\dd}{\dd t}v_i = - \frac{1}{N} \sum_{k\ne i} K_1(x_i,x_k) - \frac{1}{M} \sum_{m=1}^M K_2(x_i,d_m) - \alpha v_i, \label{ODExv}\\
 \frac{\dd}{\dd t}d_m &= u_m. \label{ODEd}
\end{align}
\end{subequations}
 In addition to the interaction terms, our model includes (linear) friction represented by the parameter $\alpha>0$. 
The individuals interact pairwise via $K_1$ with each other and with the external agents via $K_2$. The latter do not interact among themselves nor are influenced by the others. Note that such interactions can easily be considered by a simple transformation of the optimal control variables analogously to the approach in \cite{burger2009globally}. For system \eqref{longODEsys} the velocities of the external agents $u$ are assumed to be given. Later, they serve as controls for the optimization problem.


For notational convenience we further define the state vector $y := (\x,\ve,\dog) \in \mathbb{R}^{D(2N+M)}$ and the operators 
\begin{equation*}
  {\bf S}(y) = ({\bf S}_i(y))_{i=1,\dots,N},\qquad  {\bf S}_i(y) = -\frac{1}{N} \sum_{j\ne i} K_1(x_i,x_j) - \frac{1}{M} \sum_{m=1}^M K_2(x_i,d_m) - \alpha v_i. 
\end{equation*}
Using this notation, the ODE system \eqref{longODEsys} 
may be written compactly as
\begin{subequations}\label{ODEsys}
\begin{equation}\label{stateODE}
\frac{\dd}{\dd t}{y}= (\frac{\dd}{\dd t}{\x}, \frac{\dd}{\dd t}{\ve}, \frac{\dd}{\dd t}{\dog} )  = (\ve, {\bf S}(y), \con ) =: F(y,\con), 
\end{equation}
with the mapping $F: \mathbb{R}^{D(2N+M)} \times \mathbb{R}^{DM} \to \mathbb{R}^{D(2N+M)}$. The microscopic state system is supplemented with the initial conditions
\begin{equation}\label{IC_ODE}
 \x(0) = \x_0 \in \mathbb{R}^{DN}, \qquad \ve(0) = \ve_0 \in \mathbb{R}^{DN},\qquad \dog(0)=\dog_0 \in \mathbb{R}^{DM},
\end{equation}
\end{subequations}
shortly denoted  by $y(0) = y_0$.

\subsubsection{Well-posedness}
The well-posedness of \eqref{ODEsys} follows from standard theory. Here, we consider the Hilbert space  parameters, we use
\begin{align}\label{eq:space_control}
 U:=L^2(I,\mathbb{R}^{MD})
\end{align}
as the space of control, which is sufficient to prove the well-posedness of (\ref{ODExv}). Indeed, $\dog$ may be explicitly expressed as
\[
 \dog(t) = \dog_0 + \int_0^t \con(s)\dd s,
\]
which shows that $\dog$ is absolutely continuous on $I$. Consequently, we obtain the existence and uniqueness of a global solution due to the theorem of Picard and Lindel\"of and we have $\dog\in H^1(I,\mathbb{R}^{MD})$.

Altogether, we are able to define the corresponding \textit{control-to-state} operator $\mathcal{G}_N \colon U \rightarrow Y$, which maps any control parameter $\con\in U$ to the unique solution $y=\mathcal{G}_N(\con)$ of \eqref{ODEsys} in the state space $Y$ defined by
\begin{align}\label{eq:space_state_N}
 Y := H^1(I,\mathbb{R}^{ND})\times H^1(I,\mathbb{R}^{ND}) \times H^1(I,\mathbb{R}^{MD}).
\end{align}
Note, that the solution $y\in \mathcal{C}^1(I,\mathbb{R}^{ND})\times \mathcal{C}^1(I,\mathbb{R}^{ND}) \times H^1(I,\mathbb{R}^{MD})\subset Y$ is indeed in a subspace of $Y$.

\subsection{Mean-Field Model}
In order to define the limiting problem for an increasing number of individuals $N$ explicitly, we consider the empirical measure
\[
 \mu_t^N(x,v) = \frac{1}{N} \sum_{i=1}^N \delta_0(x_i(t)-x) \otimes \delta_0(v_i(t)-v).
\]
By definition, $\mu_t^N\in\mathcal{P}(\mathbb{R}^{2D})$ is a Borel probability measure that assigns the probability $\mu_t^N(A)$ of finding particles with states within a Borel measurable set $A\subset\mathbb{R}^{2D}$ in the phase space $\mathbb{R}^{2D}$ at time $t\ge 0$. If a Borel probability measure $\mu_t\in \mathcal{P}^{ac}(\mathbb{R}^{2D})$ is absolutely continuous w.r.t.~the Lebesgue measure, we denote its density by $f_t\in L^1(\mathbb{R}^{2D})$. For later use we introduce the macroscopic density $\rho_t$ of a Borel probability measure $\mu_t\in\mathcal{P}(\mathbb{R}^{2D})$ as its first marginal, i.e.,
\begin{equation*}
\rho_t(A) := \mu_t(A\times \mathbb{R}^D) = \iint_{A\times \mathbb{R}^D} f_t(x,v)\dd x\dd v,
\end{equation*}
for any Borel measurable set $A\subset\mathbb{R}^D$. The last equality holds whenever $\mu_t\in\mathcal{P}^{ac}(\mathbb{R}^{2D})$.

The link between the particle system and the mean-field equation is derived formally using ideas from \cite{Dobrushin,BraunHepp,Neunzert}. Let $h \in \mathcal{C}^\infty_c(\mathbb{R}^{2D})$ be an arbitrary smooth function on $\mathbb{R}^{2D}$ with compact support and $z_i=(x_i,v_i)\in \mathbb{R}^{2D}$. Then, it holds
\begin{equation*}
 \frac{\dd}{\dd t} \left< \mu_t^N , h \right> := \frac{\dd}{\dd  t} \frac{1}{N}\sum_{i=1}^N h(z_i) = \frac{1}{N}\sum_{i=1}^N \nabla_{x} h(z_i)\cdot\frac{\dd}{\dd t}x_i + \nabla_{v} h(z_i)\cdot\frac{\dd}{\dd t}v_i,
\end{equation*}
which allows for the formal calculation
\begin{align}\label{eq:FNweakSol}
 \langle \partial_t \mu_t^N , h \rangle 
 &= \frac{1}{N} \sum_{i=1}^N  \nabla_{x} h(z_i) \cdot v_i - \nabla_{v} h(z_i) \cdot \left( \big(K_1\ast  \rho_t^N\big)(x_i) + \frac{1}{M} \sum_{m=1}^M K_2(x_i,d_m) +\alpha v_i\right) \nonumber\\
 &= \left< \mu_t^N , \nabla_x h \cdot v - \nabla_v h \cdot \left( \big(K_1\ast  \rho_t^N\big)(x) + \frac{1}{M} \sum_{m=1}^M K_2(x,d_m) +\alpha v\right) \right>,
 \end{align}
 where 
 \[
  \big(K_1\ast  \rho_t^N\big)(x) = \int_{\mathbb{R}^{2D}} K_1(x,\bar{x}) \dd \mu_t^N(\bar{z}).
 \]
Passing to the limit $N\to\infty$ and integrating by parts, we arrive at the equation
 \[
 0 = \langle \partial_t \mu_t + v \cdot \nabla_x \mu_t + \nabla_v \cdot (S(\mu_t)  \mu_t) , h\rangle,
 \] 
 where we define
 \[
  S(\mu_t)(x,v,d) = -\big(K_1\ast  \rho_t\big)(x) - \frac{1}{M} \sum_{m=1}^M K_2(x,d_m) - \alpha v
 \]
 analogous to the particle case.
Since $h \in \mathcal{C}^\infty_c(\mathbb{R}^{2D})$ is arbitrary, we may use the variational lemma to find that
\begin{equation}\label{longPDE}
 \partial_t \mu_t  + v\cdot \nabla_x \mu_t + \nabla_v  \cdot (S(\mu_t) \mu_t )=0\qquad\text{in the sense of distributions},
\end{equation}
which is the mean-field single particle distribution, 
supplemented with the initial condition $\mu(0)=\mu_0$.

\begin{rem}
Observe that $\mu_t^N$ and $\mu_t$ satisfy exactly the same equation in the distributional sense. The above limit exists for example in the Wasserstein metric (for a detailed discussion see \cite{WassersteinConvergence}). Note, that even though $\mu_t$ has unit mass by construction, we do not need to explicitly assume that the condition is satisfied, since it is a naturally fulfilled. In particular, the control parameter does not change this property.
 
Here, we only perform the mean-field limit $N \rightarrow \infty$, while the number of external agents $M$ remains finite. 
\end{rem}

%

\subsubsection{Well-posedness} Results on the existence and uniqueness of solutions for the Vlasov equation \eqref{longPDE} can be found, e.g., in  \cite{BraunHepp,WassersteinConvergence,Dobrushin,Golse}, where the notion of solution is established in the Wasserstein space of Borel probability measures. 

\begin{defi}
 Let $\mathcal{P}_1(\mathbb{R}^{2D})$ denote the space of Borel probability measures on $\mathbb{R}^{2D}$ with finite first moment. We say that $\mu\in \mathcal{C}( I,\mathcal{P}_1(\mathbb{R}^{2D}))$ is a {\em weak measure solution} of \eqref{longPDE} with initial condition $\mu_0\in \mathcal{P}_1(\mathbb{R}^{2D})$ if for any test function $h\in \mathcal{C}_c^\infty((-\infty, T]\times \mathbb{R}^{2D})$ we have
 \[
  \int_0^T \int_{\mathbb{R}^{2D}} \big( \partial_t h_t + v\cdot\nabla_x h_t + S(\mu_t)\cdot\nabla_v h_t \big)\dd \mu_t\dd t + \int_{\mathbb{R}^{2D}} h_0\dd \mu_0 = 0.
 \]
\end{defi}

\begin{rem}
Equation \eqref{longPDE} may be equivalently expressed as a nonlinear flow
\begin{align}\label{eq:nonliner_process}
 \frac{\dd}{\dd t} Z = (\frac{\dd}{\dd t} x,\frac{\dd}{\dd t} v) = \big( v, S(\mu_t)(Z,d)\big),\qquad \mu_t=\text{law}(Z(t)),
\end{align}
 with the initial condition $Z(0)=Z_0\in\mathbb{R}^{2D}$, where $Z_0$ is a random variable distributed according to $\mu_0=\text{law}(Z_0)$. Assuming the solvability of the nonlinear process \eqref{eq:nonliner_process}, a weak measure solution of $\eqref{longPDE}$ may be represented as the push-forward of the measure along the flow $Z(t,Z_0)$, i.e., $\mu_t=Z(t,\cdot)\# \mu_0$. If we further assume that the initial measure $\mu_0\in\mathcal{P}_1(\mathbb{R}^{2D})$ has density $f_0\in L^1(\mathbb{R}^{2D})$, and that the nonlinear flow given by \eqref{eq:nonliner_process} satisfies $Z\in \mathcal{C}( I,\text{Diff}(\mathbb{R}^{2D}))$, i.e., $Z(t,\cdot)$ is a diffeomorphism for all $t\ge 0$, then $\mu_t\in\mathcal{P}_1^{ac}(\mathbb{R}^{2D})$ has density $f_t\in L^1(\mathbb{R}^{2D})$ and $f\in \mathcal{C}( I,L^1(\mathbb{R}^{2D}))$.
 \end{rem}

In order to employ the standard $L_2$-calculus we will require more regularity of the states and  assume additionally that
\begin{itemize}
 \item [\textbf{(B)}] $\mu_0\in\mathcal{P}_1(\mathbb{R}^{2D})$ has density $f_0\in L^2(\mathbb{R}^{2D})$ with compact support and the flow given by \eqref{eq:nonliner_process} satisfies $Z\in \mathcal{C}^1( I,\text{Diff}(\mathbb{R}^{2D}))$ such that $f\in H^1(I,L^2(\mathbb{R}^{2D}))$.
\end{itemize}

\begin{rem}
 Note, that the assumption on the nonlinear flow $Z\in \mathcal{C}^1( I,\text{Diff}(\mathbb{R}^{2D}))$ satisfying \eqref{eq:nonliner_process} in {\bf(B)} is not restrictive due to assumption {\bf(A)}. Indeed, $\mu\in \mathcal{C}( I,\mathcal{P}_1(\mathbb{R}^{2D}))$ implies that $S(\mu_t)$ is continuous in $t$ and locally Lipschitz in $z$. Therefore, standard ODE theory provides the required regularity for the nonlinear flow.
\end{rem}

Defining the state variable $p:= (f,\dog)$, the coupled system of \eqref{longPDE} and the ODE \eqref{ODEd} that models the movement of the external agents $\dog$ with control ${\bf w}$ may be written as
\begin{subequations}\label{PDEsys}
\begin{equation} \label{statePDE}
\partial_t p := (\partial_t f_t , \frac{\dd}{\dd t}\dog ) = -\big(\nabla_v \cdot (S(f_t)  f_t) + v \cdot \nabla_x f_t, {\bf w}\big) =: G(p,{\bf w}).
\end{equation}
The initial conditions for this system are
\begin{equation}\label{IC_PDE}
f_t|_{t=0} = f_0\in L^2(\mathbb{R}^{2D}), \qquad \dog(0) = \dog_0\in\mathbb{R}^{MD}. 
\end{equation}
\end{subequations}
In the following we shall refer to \eqref{IC_PDE} as $p(0)=p_0$. 

\begin{rem}
The problem is defined on the whole space $\mathbb{R}^{2D}$, such that  no boundary conditions are  required. In fact, we expect  from {\bf(A)} and {\bf(B)} that $f_t$  has a compact support which remains within a bounded convex domain $\Omega\subset \mathbb{R}^{2D}$ with smooth boundary for all times $t\in I$ (cf.~\cite[Thm 4.4]{FornasierSolombrino}).
\end{rem}

Finally, we define the corresponding \textit{control-to-state} operator $\mathcal{G}_\infty\colon U\to\mathcal{Y}$, $\con\mapsto p=\mathcal{G}_\infty(\con)$, where $p$ satisfies \eqref{PDEsys} and for the  state space we choose 
\begin{align}\label{eq:state_vlasov}
 \mathcal{Y}:= H^1(I,L^2(\Omega))\times H^1(I,\mathbb{R}^{MD}).
\end{align}
Here, $\Omega\subset\mathbb{R}^{2D}$ is an a priori given bounded convex domain with smooth boundary which contains the support of $f_t$ for all times $t\in I$.

\section{Control Problems and Instantaneous Control Approach}

Based on the previously introduced state systems we are now in the position to define the respective control problems. The physical task at hand is to guide the crowd of individuals such that it clusters around a certain position, while keeping its variance restricted and utilizing only limited energy for the controlling agents.

Since the question arises on the microscopic and macroscopic level, we need to have complying cost functionals modeling this issue. In fact, they need to be linked in the asymptotic limit.


%

In the following, we state one possible cost functional  using the center of mass and the variance of the crowd, which meets these requirements. Further, we are  describing the instantaneous control procedure that we are going to apply to control the crowd.

\subsection{Cost Functional in the Microscopic Setting}
Again, let  $\x$ denote the particle positions. Then,  the center of mass $\mathbb{E}_N$ and variance $\mathbb{V}_N$ for the particles are given by
\begin{equation*}
 \mathbb{E}_N(\x(t)) = \frac{1}{N}  \sum_{i=1}^N x_i(t),\qquad \mathbb{V}_N(\x(t)) = \frac{1}{N} \sum_{i=1}^N |x_i(t)-\mathbb{E}_N(\x(t))|^2,
\end{equation*}
and we define the cost functional
\begin{equation}\label{eq:cost functional_N}
 J_N(y(\con),\con) = \int_0^T \frac{\sigma_1}{4T} |\mathbb{V}_N(\x(t)) - \bar{V}_N|^2+ \frac{\sigma_2}{2T}  |\mathbb{E}_N(\x(t))-\Edes|^2  + \frac{\sigma_3}{2M} \norm{\con(t)}_{\mathbb{R}^{MD}}^2 \dd t,
\end{equation}
where $\bar{V}_N >0$ is a given desired variance value. 

The first part of the functional penalizes the spread of the particles from the desired variance $\bar{V}_N$, while the second term measures the distance of the center of mass from  the desired destination $\Edes\in\mathbb{R}^D$. The third term measures the control costs in terms of the kinetic energy of the agents. Mathematically, it introduces more convexity in the cost functional. The positive parameters $\sigma_i$ allow to adjust the influence of the different parts of the cost functional. 




\subsection{Cost Functional in the Mean-Field Setting}
As already discussed, the control problem on the microscopic level has to match the one on the mean-field level as $N \rightarrow \infty$. Thus, we choose the following cost functional
\begin{equation}\label{eq:cost functional_MF}
 J_\infty(p({\bf w}),{\bf w}) =  \int_0^T \frac{\sigma_1}{4T} |\mathbb{V}_{\infty}(f_t) - \bar{V}_\infty |^2+ \frac{\sigma_2}{2T}  |\mathbb{E}_{\infty}(f_t)-\Edes|^2 + \frac{\sigma_3}{2M} \norm{\textbf{w}_t}_{\mathbb{R}^{MD}}^2 \dd t,
 \end{equation}
 in which  the center of mass $\mathbb{E}_\infty$ as well as the  variance $\mathbb{V}_\infty$ only depend on the macroscopic density and are defined as 
\begin{equation*}
 \mathbb{E}_{\infty}(f_{t}) = \iint_{\mathbb{R}^{D}\times \mathbb{R}^D} x  f_t\dd x \dd v,\qquad \mathbb{V}_{\infty}(f_{t}) = \iint_{\mathbb{R}^{D} \times \mathbb{R}^D} |x -\mathbb{E}_{\infty}(f_{t})|^2 f_t\dd x \dd v.
\end{equation*}

Note that the microscopic cost functional can be derived from the mean-field cost functional by using the empirical measure $\mu^N$. In fact, it holds
\begin{gather*}
 \mathbb{E}_{\infty}(f^N_{t}) = \iint_{\mathbb{R}^{D}\times \mathbb{R}^D} x  f^N_t\dd x \dd v = \frac{1}{N}\sum_{i=1}^N x_i(t), \\ 
\intertext{and}  
 \mathbb{V}_{\infty}(f^N_{t}) = \iint_{\mathbb{R}^{D} \times \mathbb{R}^D} |x -\mathbb{E}_{\infty}(f^N_{t})|^2 f^N_t\dd x \dd v = \frac{1}{N} \sum_{i=1}^N(x_i(t) - \mathbb{E}_N(\x(t)))^2.
\end{gather*}

\subsection{Instantaneous Control}
The idea of the instantaneous control approach is to split the time interval of interest into several smaller time intervals and to solve the control problem sequentially on these subintervals. Clearly this reduces the memory consumption of the algorithm significantly.
Moreover, it is in our setting more realistic than optimal control over the whole time horizon, because we expect the agents to forecast the movement of the crowd only on short time intervals. 

Hence, we split time interval of investigation into $K+1$ subintervals $I_k=[t_{k-1},t_k]$, where $t_0 = 0$ and $t_K = T$, $k=1,\dots,K$ such that $[0,T] = \bigcup_{k=1}^K[t_{k-1},t_k]$.  The corresponding spaces $Y^k, \mathcal{Y}^k$ and $\mathcal{U}_{ad}^k$ are precisely defined  below. Then, we study the family of optimization problems on each time slice $I_k$ on the microscopic level given by
\vspace{0.3cm}
\begin{mdframed}
\begin{problem}\label{Opt_ODE}
 For each subinterval $I_k$, $k=1,\dots,K$, find an optimal pair $(y^*_k,\con^*_k) \in Y\times \mathcal{U}_{ad}$ such that
 \begin{gather*} 
  (y_k^*,\con_k^*) = \argmin\nolimits_{(y,\con) \in Y^k \times \mathcal{U}_{ad}^k}  J_N(y,\con) \\ \text{ subject to } \\ \text{ IVP } \eqref{ODEsys} \text{ on } I_k  \\
  \text{and initial data $y_k(t_{k-1}) = y_{k-1}^*(t_{k-1})$, where} \\ 
  J_N(y,\con) =  \int_{I_k} \frac{\sigma_1}{4T} |\mathbb{V}_N(\x(t)) - \bar{V}_N|^2+ \frac{\sigma_2}{2T}  |\mathbb{E}_N(\x(t))-\Edes|^2  + \frac{\sigma_3}{2M} \norm{\con(t)}_{\mathbb{R}^{MD}}^2 \dd t.
 \end{gather*}
\end{problem}
\end{mdframed}
\vspace{0.3cm}
On the mean-field level we obtain the following optimization problem:
\vspace{0.3cm}
\begin{mdframed}
\begin{problem}\label{Opt_PDE}
 For each subinterval $I_k$, $k=1,\dots,K$, find an optimal pair $(p_k^*,{\bf w}_k^*) \in \mathcal{Y}^k\times\mathcal{U}_{ad}^k$ such that
 \begin{gather*} 
  (p_k^*,{\bf w}_k^*) = \argmin\nolimits_{(p,{\bf w}) \in \mathcal{Y}^k\times \mathcal{U}_{ad}^k}  J_\infty(p,{\bf w})\\ \text{ subject to} \\  \text{system \eqref{PDEsys} on } I_k \text{ and initial data $p_k(t_{k-1}) = p_{k-1}^*(t_{k-1}),$ }, \text{where} \\
  J_\infty(p,{\bf w}) = \int_{I_k} \frac{\sigma_1}{4T} |\mathbb{V}_{\infty}(f_t) - \bar{V}_\infty |^2+ \frac{\sigma_2}{2T}  |\mathbb{E}_{\infty}(f_t)-\Edes|^2 + \frac{\sigma_3}{2M} \norm{\textbf{w}_t}_{\mathbb{R}^{MD}}^2 \dd t.
 \end{gather*}
 \end{problem}
 \end{mdframed}

\begin{rem}
For the first iterate the initial data is given by $y_1^*(0) = y_0$ and   $p_1^*(0) = p_0$, respectively.
\end{rem}

Further, we restrict the velocities of the external agents and  define the spaces of admissible controls $\mathcal{U}_{ad}^k$, for $k=1,\dots,K$, as 
\begin{equation}\label{eq:control}
 \mathcal{U}_{ad}^k := \left\{ \con \in L^2(I_k,\mathbb{R}^{MD})\; \colon\; |u_m(t)| \le u_{\text{max}}, \quad m=1,\dots,M \right\}.
\end{equation}

\subsubsection{Existence of controls}
An existence result for the optimization problem on each time slice may be deduced in a straight-forward way from the results in \cite{FornasierSolombrino}. Even though the results are based on a sparse control setting, similar arguments may be applied in the present setting as well. 

\begin{thm}
Assume {\bf (A)} and {\bf (B)}. Then, there exists on each time interval $I_K$, $k =1,\ldots,K$, an optimal pair $(y^*_k,\con^*_k) \in Y\times \mathcal{U}_{ad}$ for Problem 1 and an optimal pair $(p_k^*,{\bf w}_k^*) \in \mathcal{Y}^k\times\mathcal{U}_{ad}^k$ for Problem 2, respectively.
\end{thm}

\begin{rem}
The choice of initial conditions for each subinterval $I_k$ indicates that we obtain a sub-optimal solution for the problem on the interval $[0,T]$ by cluing the optimal solutions of the subintervals together.
\end{rem}

\section{First Order Necessary Conditions} \label{sec:adjoint}
We are going to apply adjoint based descent methods to solve the control problems on each time slice. Therefore, we formally derive the adjoints and the optimality conditions for the particle and the mean-field optimal control problem with the  help of the extended Lagrange functional. Furthermore, we introduce the reduced cost functional and its gradient. Throughout this section we consider an arbitrary time slice $I$ and denote its right-end point $I_R.$

\subsection{Adjoint of the microscopic problem}
Let the control space $U$ and state space $Y$ be the Hilbert spaces
\begin{equation*} 
 U=L^2( I_k,\mathbb{R}^{MD}),\qquad Y = [H^1( I,\mathbb{R}^{ND})]^2 \times H^1( I,\mathbb{R}^{MD}),
\end{equation*}
with $\mathcal{U}_{ad}\subset U$ defined in \eqref{eq:control}. Further, we define $X:=[L^2( I,\mathbb{R}^{ND})]^2\times L^2( I,\mathbb{R}^{MD})$ and
\[
 Z:=X\times\big([\mathbb{R}^{ND}]^2\times \mathbb{R}^{MD}\big),
\]
as the space of Lagrange multipliers. 

The state operator $e_N\colon Y \times U \rightarrow Z^*$ of the microscopic problem is given by
\begin{equation*}
e_N(y,\con) = \begin{pmatrix} \frac{\dd}{\dd t}y - F(y,\con) \\ y(0)-y_0 \end{pmatrix}
\end{equation*}
or in weak form
\begin{equation*}
 \langle e_N(y,\con),(\xi,\eta) \rangle_{Z^*,Z} = \int_{I} (\frac{\dd}{\dd t}y(t) - F(y(t),\con(t))) \cdot \xi(t) \dd t + (y(0) - y_0) \cdot \eta.
\end{equation*}
Let $(\xi,\eta)\in Z$ denote the Lagrange multiplier, which is in fact the adjoint state. Then, the extended Lagrange functional corresponding to Problem \ref{Opt_ODE} reads
\begin{equation*}
 \mathcal{L}_{N}(y,\con,\xi,\eta)  = J_{N}(y,\con) + \langle e_N(y,\con),(\xi,\eta) \rangle_{Z^*,Z}.
\end{equation*}
As usual the first-order optimality condition of Problem \ref{Opt_ODE} is derived by solving
\begin{equation*}
 d \mathcal{L}_N(y,\con,\xi,\eta) =0. 
\end{equation*}
The derivative w.r.t.~the adjoint state results in the state equation, while the derivative with respect to the state $y$ yields the adjoint system and the optimality condition is obtained by the derivative w.r.t.~the control (for details see, e.g., $\con$ \cite{Pinnau}). 

For the calculations we denote the three parts of the cost functional by $J^i$, $i=1,2,3$, with
\begin{gather*}
 J_N^1(y) = \frac{\sigma_1}{4T} \int_{I} |\mathbb{V}_N(\x(t)) - \bar{V}_N|^2 \dd t, \qquad J_N^2(y) =  \frac{\sigma_2}{2T} \int_{I} \norm{\mathbb{E}_N(\x(t))-\Edes}_{\mathbb{R}^{D}}^2 \dd t,\\
  J_N^3(\con)=\frac{\sigma_3}{2M} \int_{I} \norm{\con(t)}_{\mathbb{R}^{MD}}^2 \dd t.
\end{gather*}

For any $h=(h^y=(h^x,h^v,h^d),h^u)\in Y\times U$, the following G\^ateaux derivatives are obtained
\begin{align*}
 d_{\x} J_N^1(y) [h^x] 
 &=  \frac{\sigma_1 }{NT} \int_{I}  (\mathbb{V}(\x(t)) - \bar{V}_N)(\x(t) - \mathbb{E}(\x(t))) \cdot h^x(t) \dd t,  \\
 d_{\x} J_N^2(y) [h^x] 
 &= \frac{\sigma_2}{NT} \int_{I} (\mathbb{E}(\x(t)) - \Edes) \cdot  h^x(t) \dd t, \\
 d_{\con} J_N^3(\con) [h^u] &=  \frac{\sigma_3}{M} \int_{I} \con(t) \cdot h^u(t) \dd t.
\end{align*}
For the second part of the Lagrangian we derive
\begin{align*}
 \langle d_\x e_N(y,\con)[h^x],(\xi,\eta) \rangle &= \int_{I} \frac{\dd}{\dd t}h^x(t)\cdot\xi_1(t) - d_{\x} {\bf S}(y)[h^x(t)]\cdot \xi_2(t) \dd t + h^x(0)\cdot\eta_1, \\
 \langle d_{\ve} e_N(y,\con)[h^v], (\xi,\eta) \rangle &= \int_{I} \left(\frac{\dd}{\dd t}h^v(t) + \alpha h^v(t)\right)\cdot \xi_2(t) - h^v(t) \cdot \xi_1(t)  \dd t + h^v(0)\cdot\eta_2, \\
 \langle d_{\dog} e_N(y,\con)[h^d], (\xi,\eta) \rangle &= \int_{I} \frac{\dd}{\dd t}h^d(t)\cdot\xi_3(t) - d_{\dog}{\bf S}(y)[h^d(t)] \cdot \xi_2(t) \dd t + h^d(0)\cdot\eta_3, \\
 \langle d_{\con} e_N(y,\con)[h^u] , (\xi,\eta) \rangle &= -\int_{I} h^u(t)\cdot \xi_3(t) \dd t.
\end{align*}

 Assuming that $\xi\in Y$ one may formally derive the strong formulation of the adjoint system. Indeed, integrating by parts and using the fact that $d_{\x} {\bf S}(y)$ and $d_{\dog}{\bf S}(y)$ are symmetric matrices, we arrive at the following result.

\begin{prop}\label{ODEKKT}
Let $(y_*,u_*)$ be an optimal pair for Problem 1. Then,  the first-order optimality condition corresponding to Problem \ref{Opt_ODE} reads
 \begin{equation}\label{eq:var_ode}
  \int_{I} \left(\frac{\sigma_3}{M} \con_*(t) - \xi_3(t) \right)\cdot(\con(t)-\con_*(t)) \dd t\ge 0  \qquad \text{for all\; $\con\in \mathcal{U}_{ad}$},
 \end{equation}
 where $\xi=(\xi_1,\xi_2,\xi_3)\in Y$ satisfies the adjoint system given by
 \begin{subequations}\label{ad_ODE}
 \begin{equation}\label{xi}
  \frac{\dd}{\dd t}\xi_1 = -d_{{\emph \x}} {\bf S}(y_*)[\xi_2]-d_{\emph \x}J_N(t), \qquad \frac{\dd}{\dd t}\xi_2 = \xi_1 - \alpha \xi_2, \qquad  \frac{\dd}{\dd t}\xi_3 = -d_{{\emph \dog}} {\bf S}(y_*)[\xi_2],
 \end{equation}
 with
  \begin{equation}\label{terminal_ODE}
  d_{\emph \x} J_N(t) = \frac{\sigma_1}{NT}\Big((\mathbb{V}({\emph \x_*}(t)) - \bar{V}_N)({\emph \x_*}(t) - \mathbb{E}_N({\emph \x_*}(t)))\Big) + \frac{\sigma_2}{NT} \Big(\mathbb{E}_N({\emph \x_*}(t))-\Edes\Big),
 \end{equation}
supplemented with the terminal conditions $\xi_1(I_R) = 0$, $\xi_2 (I_R) = 0$, $\xi_3(I_R) = 0$.

 \end{subequations}
\end{prop}

\subsection{Adjoint of the Mean-Field Problem}
 Here, we assume that $p = (f,\dog)$ lies within the state space $\mathcal{Y}$ of the PDE optimization problem, where
 \[
  \mathcal{Y} = H^1( I,L^2(\Omega)) \times H^1( I,\mathbb{R}^{MD}).
 \]
 Let $\mathcal{X}:=H^1( I,L^2(\Omega))\cap L^2( I,H^1(\Omega))\times L^2( I,\mathbb{R}^{MD})$ and set $\mathcal{Z} := \mathcal{X}\times \left( L^2(\Omega) \times \mathbb{R}^{MD} \right)$ to be the space of adjoint states with dual $\mathcal{Z}^*$. Note that the control space $U$ is identical to the one of the particle problem. 
 
 Now, we define the mapping $e_\infty\colon \mathcal{Y} \times U \rightarrow \mathcal{Z}^*$ by
\begin{align*}
 \langle e_\infty(p,{\bf w}),(\varphi,\eta) \rangle_{\mathcal{Z}^*, \mathcal{Z}} &= -\int_{I} \int_\Omega \big( \partial_t g_t + v\cdot\nabla_x g_t + S(f_t)\cdot\nabla_v g_t \big)f_t \dd z \dd t + \int_{I} (\frac{\dd}{\dd t}\dog - {\bf w})\cdot \varphi_d \dd t \\
 &\hspace*{2em}+ \int_\Omega g(T)f(T) - g(0)f(0) \dd z - \int_\Omega (f(0) -f_0) \eta_f \dd z   + (\dog(0)-\dog_0)\cdot \eta_d,
\end{align*}
with the adjoint state $(\varphi,\eta)\in \mathcal{X}\times \left( L^2(\Omega) \times \mathbb{R}^{MD} \right)$, $\varphi=(g,\varphi_d)$ and $\eta=(\eta_f,\eta_d)$.

 Similar to the microscopic case we define the extended Lagrangian corresponding to Problem \ref{Opt_PDE} as
\begin{equation*}
 \mathcal{L}_\infty(p,{\bf w},\varphi,\eta) = J_\infty(p,{\bf w}) + \langle e_\infty(p,{\bf w}),(\varphi,\eta) \rangle_{\mathcal{Z}^*, \mathcal{Z}},
\end{equation*}

Again, we introduce the three parts of the cost functional by $J^i$, $i=1,2,3$, via
\begin{gather*}
 J_\infty^1(p) = \frac{\sigma_1}{4T} \int_{I} |\mathbb{V}_\infty(f_t) - \bar{V}_\infty|^2 \dd t, \qquad J_\infty^2(p) =  \frac{\sigma_2}{2T} \int_{I} |\mathbb{E}_\infty(f_t)-\Edes|^2 \dd t, \\
  J_\infty^3({\bf w})=\frac{\sigma_3}{2M} \int_{I} \norm{\textbf{w}(t)}_{\mathbb{R}^{MD}}^2 \dd t.
\end{gather*}

Analogously to the microscopic case we derive the adjoint system and the optimality condition by calculating the derivatives of $\mathcal{L}_\infty$ w.r.t.~the state variable $p$ in direction $h^p=(h^f,h^d)\in\mathcal{Y}$, and the control $w$ in direction $h^w\in U$. The standard $L_2$-calculus yields
\begin{align*}
 d_f J_\infty^1 (p) [h^f] 
 &=\frac{\sigma_1}{T} \int_{I} \int_\Omega \Big(\mathbb{V}_\infty(f_t) - \bar{V}_\infty\Big) |x - \mathbb{E}_\infty(f_t)|^2  \dd h_t^f \dd t, \\
 d_f J_\infty^2 (p) [h^f] 
 &= \frac{\sigma_2}{T} \int_{I} \int_\Omega x\cdot (\mathbb{E}_\infty(f_t) - \Edes) \dd h_t^f \dd t, \\
 d_{{\bf w}} J_\infty^3 ({\bf w}) [h^w] &= \frac{\sigma_3}{M} \int_{I} {\bf w}(t) \cdot h^w(t) \dd t.
\end{align*}
Let $\varphi = (g, \varphi_d)$ be the adjoint state corresponding to $p=(f,\dog)$. Then, we obtain for the constraint part in the extended Lagrange functional the following G\^ateaux derivatives:
\begin{align*}
 \langle d_f e_\infty(p,{\bf w})[h^f] , (\varphi,\eta) \rangle  &=  -\int_{I} \int_\Omega  \big( \partial_t g_t + v\cdot\nabla_x g_t + S(f_t)\cdot\nabla_v g_t \big)h_t^f \dd z \dd t, \\
 &\hspace*{4em}- \int_{I} \int_\Omega d_fS(f_t)[h_t^f]\cdot\nabla_vg_t \, f_t\dd z \dd t \\ 
 &\hspace*{4em}+ \int_\Omega g(T) h^f(T) -  h^f(0) g(0) - h^f(0) \eta_f \dd z,\\
 \langle d_{\dog} e_\infty(p,{\bf w})[h^d] , (\varphi,\eta)  \rangle &= \int_{I} \frac{\dd}{\dd t}h^d(t)\cdot \varphi_d(t) \dd t + h^d(0)\cdot \eta_d, \\
 &\hspace*{4em}- \int_{I} \int_\Omega d_{\dog} S(f_t)[h^d(t)]\cdot \nabla_vg_t \, f_t\dd z \dd t, \\
 \langle d_{{\bf w}}  e_\infty(p,{\bf w})[h^w] , (\varphi,\eta)  \rangle  &= - \int_{I} h^w(t) \cdot \varphi_d(t) \dd t.
\end{align*}
%
Assuming again that the adjoint state $\varphi_d$ is sufficiently regular, we may integrate by parts to obtain a strong formulation of the adjoint system. For the terms involving derivatives of $S(f_t)$ we calculate  the following representations:
\begin{align}
 \int_\Omega d_fS(f_t)[h_t^f]\cdot\nabla_vg_t \, f_t\dd z
 &= -\int_\Omega \int_\Omega K_1(x,\bar x)\, h_t^f \dd \bar z \cdot \nabla_v g_t(z)\,f_t \dd z  \nonumber \\
 &= \int_\Omega \left(\int_\Omega K_1(\bar x,x)  \cdot \nabla_v g_t(z)\,f_t \dd z\right) h_t^f \dd \bar z \nonumber \\
 &=: \int_\Omega D_f(f_t)[g_t](z)\,h_t^f\dd z, \nonumber \\
 \int_\Omega d_{\dog} S(f_t)[h^d(t)]\cdot \nabla_vg_t \, f_t\dd z &= \int_\Omega \left( d_\dog S(f_t)[h^d(t)]\right) \cdot \nabla_v g_t(z)\,f_t \dd z \nonumber \\
 &= \left(\int_\Omega d_\dog S(f_t)[\nabla_vg_t(z)]\,f_t \dd z\right)\cdot h^d(t)\nonumber \\
 \label{D_fS}
 &=: D_\dog (f_t)[g_t]\cdot h^d(t). 
\end{align}
This yields the following adjoint system and optimality condition.

\begin{prop}\label{PDEKKT}
Let $(p_*,w_*)$ be an optimal pair for Problem 2. Then, the optimality condition corresponding to Problem \ref{Opt_PDE} reads
 \begin{equation}\label{eq:var_vlasov}
  \int_{I} \left(\frac{\sigma_3}{M} {\bf w}_*(t) - \varphi_d(t) \right)\cdot({\bf w}(t)-{\bf w}_*(t)) \dd t\ge 0  \qquad \text{for all\; ${\bf w} \in \mathcal{U}_{ad}$},
 \end{equation}
 where $\varphi=(g,\varphi_d)\in \mathcal{Y}$ satisfies the adjoint system given by
 \begin{subequations}\label{adSys}
 \begin{align}
  \partial_t g_t + v\cdot \nabla_x g_t &= -S(f_{t,*}) \cdot \nabla_v g_t + D_f(f_{t,*})[g_t] - d_{\emph \x} J_\infty(t)  \label{ad_PDE}, \\
  \frac{\dd}{\dd t}\varphi_d &= -D_{\emph\dog}(f_{t,*}) [g_t], \label{ad_dog}
 \end{align}
 where
 \begin{equation}\label{dxJinf}
 d_{\emph \x} J_\infty(t) = \frac{\sigma_1}{T} \Big( \mathbb{V}(f_{t,*}) - \bar{V}_\infty \Big) |x-\mathbb{E}_\infty(f_{t,*})|^2 + \frac{\sigma_2}{T} (\mathbb{E}_\infty(f_{t,*}) - \Edes) \cdot x \qquad\text{on\; $\Omega$}.
\end{equation}
 The system is supplemented with the terminal conditions $g_{I_R} = 0$ and $\varphi_d({I_R}) =  0$.

\end{subequations}
\end{prop}

\begin{rem}
 Note, that the optimality conditions of Problem \ref{Opt_ODE} and \ref{Opt_PDE} coincide. This is due to the fact that $J^3$ and the set of admissible controls are identical in both problems. The mean-field limit only affects the adjoint system, but formally they can be identified in the following way: along the characteristics of the partial differential equation for $g$, $\xi_1$ corresponds to $\frac{1}{N}\nabla_x g$ and $\xi_2$ to $\frac{1}{N}\nabla_v g$. Finally, $\varphi_d$ can be identified directly with $\xi_3$. \cite{Schafe2}
\end{rem}

\begin{rem}
 The variational inequalities \eqref{eq:var_ode} and \eqref{eq:var_vlasov} may be equivalently expressed as fixed point problems in terms of a projection operator $\text{Proj}_U\colon U\to \mathcal{U}_{ad}$ which is defined by (see  \cite{Pinnau})
 \[
  \text{Proj}_U(h)=\argmin\nolimits_{u \in \mathcal{U}_{ad}} \|u-h\|_U\qquad\text{for any\; $h\in U$}.
 \]
 Consequently, the variational inequalities \eqref{eq:var_ode} and \eqref{eq:var_vlasov} may be expressed as
 \begin{align*}
  u_* = \text{Proj}_U(u_*-\gamma k)\in\mathcal{U}_{ad},
 \end{align*}
 where $k(\con)=\sigma_3/M \con_* - \xi_3$ for the microscopic case \eqref{eq:var_ode} and $k({\bf w})=\sigma_3/M {\bf w}_* - \varphi_d$ for the mean-field case \eqref{eq:var_vlasov}. 
 
 In our particular case, $\text{Proj}_U$ has the explicit representation given by
 \begin{equation}\label{projection}
  \text{Proj}_U(h)(t) = \begin{cases}
  u_{\text{max}}\frac{h_m(t)}{|h_m(t)|} & \text{for\; $|h_m(t)|>u_{\text{max}}$}, \\
  h_m(t) & \text{otherwise},
  \end{cases}\quad m=1,\ldots,M,\quad \text{a.e. in\; $ I$}.
 \end{equation}
\end{rem}

Due to the high dimensionality of the problem, it is not recommended solving the system consisting of state equation, adjoint equation and optimality condition all at once. Therefore, we employ a gradient descent method. In order to define this method we employ the gradient of the reduced cost functional which is derived in the following for both the microscopic problem and its mean-field counterpart.

\subsection{Gradient of the Reduced Cost Functional} \label{sec:redGrad}
In this section we introduce the reduced cost functionals $\hat{J}_N(\con)$ and $\hat{J}_\infty({\bf w})$ and formally calculate their gradients $\nabla \hat{J}_N(\con)$ and $\nabla \hat{J}_\infty({\bf w})$ which we need for the descent algorithms. 

By means of the control-to-state operators $\mathcal{G}_N\colon U\to Y$ and $\mathcal{G}_\infty\colon U\to\mathcal{Y}$ introduced in Section~\ref{sec:state}, we define the reduced cost functionals as
\begin{align*}
 \hat{J}_N(\con) := J_N(\mathcal{G}_N(\con),\con),\qquad
 \hat{J}_\infty({\bf w}) := J_\infty(\mathcal{G}_\infty({\bf w}),{\bf w}).
\end{align*}
Assuming sufficient regularity for $\mathcal{G}^N$ and $\mathcal{G}^\infty$ we further derive the gradients of the reduced cost functionals. Making use of the state equation $e_N(y,\con) = 0$  and $e_\infty(p,{\bf w})=0$ we implicitly obtain $d\mathcal{G}_N(\con)$ and $d\mathcal{G}_\infty({\bf w})$ via
\begin{align*}
0= d_\con e_N(\mathcal{G}_N(\con),\con) &= d_y e(\mathcal{G}_N(\con),\con)[d\mathcal{G}_N(\con)] + d_\con e_N(\mathcal{G}_N(\con),\con), \\ 
0= d_{\bf w} e_\infty(\mathcal{G}_\infty({\bf w}),{\bf w}) &= d_p e(\mathcal{G}_\infty({\bf w}),{\bf w})[d\mathcal{G}_\infty({\bf w})] + d_{\bf w} e_\infty(\mathcal{G}_\infty({\bf w}),{\bf w}),
\end{align*}
With the help of the adjoint equations
\begin{equation*}
 (d_ye(y,\con))^*[\xi] = - d_y J_N(y,\con)   \qquad \text{and} \qquad (d_pe(p,{\bf w}))^*[\varphi] = - d_p J_\infty(p,{\bf w})
\end{equation*}
we may calculate the G\^ateaux derivative of $\hat{J}$ in the direction $h\in U$, which gives
\begin{align*}
 d\hat{J}_N(\con)[h] &= \langle d_y J_N(y,\con), d\mathcal{G}_N(\con)[h]\rangle_{Y^*,Y} + \langle d_\con J_N(y,\con),h\rangle_U = \langle \frac{\sigma_3}{M} \con -\xi_3, h \rangle_{U},\\
 d\hat{J}_\infty({\bf w})[h] &= \langle d_p J_\infty(p,{\bf w}), d\mathcal{G}_\infty({\bf w})[h]\rangle_{\mathcal{Y}^*,\mathcal{Y}} + \langle d_{\bf w} J_\infty(p,{\bf w}),h\rangle_U =\langle \frac{\sigma_3}{M} {\bf w} -\varphi_d,h  \rangle_{U}.
\end{align*}
Since $U$ is a Hilbert space, we may use the Riesz representation theorem to identify the gradients
\begin{equation}\label{gradient}
 \nabla \hat{J}_N (\con) = \frac{\sigma_3}{M} \con -\xi_3 \qquad \text{and} \qquad \nabla \hat{J}_\infty ({\bf w}) = \frac{\sigma_3}{M} {\bf w} - \varphi_d.
\end{equation}
\section{Numerical Schemes and Instantaneous Control}\label{sec:Alg}

For the numerical solution of the two control problems we use the instantaneous control algorithm.  In this case the time slices are equally sized and the size coincides with the time step size of the numerical scheme, i.e., the gradient is based on the information of one time step only. Further, we discuss the numerical schemes for the solution of the respective forward and backward differential equations.


\subsection{Numerics for the Forward and Adjoint IVP}
We compute the state variable with one step of the explicit fourth order Runge--Kutta solver for the IVP \eqref{ODEsys}. For the mean-field limit we need a numerical scheme which is independent of the number of particles involved. Hence, we rescale the adjoint ODE by multiplying with $N$. This has the effect that the $N$-dependence of the terms in \eqref{terminal_ODE} emerging from the cost functional vanishes. 

In fact, for $i=1,\dots,N$, we set $r_i(t) = N \xi_i^1(t)$ and $s_i(t) = N \xi_i^2(t)$, and obtain the rescaled adjoint ODE system
\begin{subequations}\label{rescaled_ad}
\begin{align}
 \frac{\dd}{\dd t}r_i &= - \frac{1}{N} \sum_{j\ne i} \nabla _{x_i} K_1(x_i,x_j)(s_i-s_j) - \frac{1}{M} \sum_m\nabla_{x_i}K_2(x_i,d_m) s_i - \frac{1}{T}d_{x_i}J_N(t), \\
 \frac{\dd}{\dd t}s_i &= - r_i- \alpha s_i, \\
 \frac{\dd}{\dd t}\varphi_i &= \frac{1}{NM} \sum_{i=1}^N\nabla_{x_i} K_2(x_i,d_m) s_i,
\end{align}
where
\begin{equation}
d_{x_i} J_N(t) = \sigma_1(\mathbb{V}(\x(t)) - \bar{V}_N)\big(x_i(t) - \mathbb{E}_N({\x}(t))\big) + \sigma_2\big(\mathbb{E}_N({ \x}(t))-\Edes\big),
\end{equation}
with terminal conditions $r(T) = 0$, $s(T)=0$ and $\varphi(T)=0$.
\end{subequations}
We apply one step of the explicit fourth order Runge-Kutta scheme to \eqref{rescaled_ad} for computing the adjoint. 

Note, that the main computational effort for a large particle number $N$ comes from the interaction potentials.

\subsection{Numerics for the Mean-field Equation and its Adjoint}
The forward and backward solves for the mean-field optimization are realized using a {\em Strang splitting} scheme\cite{StrangSplitting}. This scheme applies a semi-Lagrangian solver in the spatial direction and a semi-implicit finite-volume scheme in the velocity space. Using the following short hand notation for  \eqref{PDEsys}
\begin{equation*}
\partial_t f = -v \cdot \nabla_x f - \nabla_v  \cdot (S(f)  f),
\end{equation*}
we define the splitting
\begin{subequations}\label{discretization}
\begin{align}\label{dis_1}
\partial_t f^* &= -\frac{1}{2} \nabla_v \cdot ( S(f^*)   f^*) , &f^*(t) &= f(t), \\
\label{dis_2} \partial_t f^{**} &= - v \cdot \nabla_x f^{**}, &f^{**}(t) &=f^*(t+\tau), \\
\label{dis_3} \partial_t f &= -\frac{1}{2} \nabla_v  \cdot (S(f) f) , &f(t) &= f^{**}(t+\tau).
\end{align}
\end{subequations}
A Semi-Lagrangian method \cite{sonnendrucker1999semi,klar2009semi} is used to solve \eqref{dis_2}, that means the computations rely on a fixed grid and yet we make use of the Lagrangian ansatz and consider the transport along characteristics.  To obtain the characteristic curves we solve ODEs using a second order Runge-Kutta scheme. The inital point for each transport step is a grid point. Since we cannot ensure that also the endpoint of each transport step is again grid point, we need to interpolate the data to the grid. This interpolation is realized by a polynomial reconstruction based on cubic Bezier curves, which is again of second order. 

For the discretization in velocity space, \eqref{dis_1} and \eqref{dis_3}, we adopt a second order finite volume scheme where the advection is approximated by a Lax-Wendroff flux \cite{leveque2002finite,quarteroni2008numerical}. Oscillations caused by non-smooth solutions are intercepted with the help of a van-Leer limiter \cite{VanLeerLimiter}. More details on this second order scheme can be found in \cite{Andreas}. 

Basically the same code is used for the adjoint system: we rewrite
\begin{equation*}
S(f) \cdot \nabla_v g = -\nabla_v \cdot (S(f) g) + 2 \alpha g.
\end{equation*}
Further, we need to add the term resulting from the linearization of the non-linear interaction and the cost functional. Altogether we obtain the splitting for the adjoint system
\begin{align*}
\partial_t g^* &= -d_x J_\infty + \frac{1}{2}  \left( -\nabla_v \cdot (S(f) g^*) + 2 \alpha g^* + D_f(f_t)[g_t]\right) , &g^*(t) &= g(t), \\
 \partial_t g^{**} &= - v \cdot \nabla_x g^{**}, &g^{**}(t) &=g^*(t+\tau), \\
\partial_t g &= \frac{1}{2} \left( -\nabla_v \cdot (S(f) g) + 2 \alpha g  + D_f(f_t)[g_t] \right), &g(t) &= g^{**}(t+\tau).
\end{align*}
with $D_f(f_t)[g_t]$ as defined in \eqref{D_fS} and $d_x J_\infty$ as in \eqref{dxJinf}.

\begin{rem}
Note, that the forward as well as the backward system are high-dimensional, since for our example in $\R^2$ we need in fact to solve equations in $\R^4$.  This yields an immense computational effort and is also challenging for the data handling. The last is the essential bottleneck for an optimal control approach, since then one needs to store the whole time evolution of the forward problem to do one backward solve for the evaluation of the gradient. Nevertheless, numerical results in the optimal control setting on a short time horizon can be found in \cite{Diss}. 
\end{rem}

\subsection{Instantaneous Control Algorithm}
We use the steepest descent steps to update the control once on every time slice, i.e.,
\begin{equation}\label{eq:update}
 \tilde{c}_{k+1} = c_k - \omega_k q_k,
\end{equation} 
where $c_k$ denotes the current control and  $q_k = \nabla \hat{J}(u_k)$. 
Due to the uniform bound on the velocity  we need to ensure that the  computed the controls are feasible. Therefore, we  project the controls onto the feasible set in every iteration using the operator $\text{Proj}_U$ defined in \eqref{projection}. The step sizes $\omega_k$ are obtained with help of the projected Armijo step size rule \cite{Pinnau} (see Algorithm \ref{ArmRule}).
\RestyleAlgo{boxruled}
\begin{algorithm2e}[h!]
\caption{Projected Armijo Stepsize Rule }\label{ArmRule}
 \KwData{Current control $c_k$, gradient $q_k$, initial $\omega_0$, initial $\gamma$}
 \KwResult{new control $c_{k+1}$}
 initialization\;
 \While{ $\hat{J}( \emph{Proj}_U[c_k - \omega_k q_k]) \ge  \hat{J}(c_k) - \gamma\omega_k \norm{q_k}^2$ }{
   $\omega_k = \omega_k/2$
 }
\end{algorithm2e}

The control on time slice $k+1$ is initialized with
\begin{equation}\label{IC_update_control}
	c_{k+1} = 0.1 c_k,
\end{equation}
which ensures  that the cost term of the cost functional is initially non-zero. Hence, changing the direction of the controls while fixing the speed does not increase the third part of the cost functional. Therefore, it is easier to find descent directions. 

We combine the solvers for the state system, the adjoint system and the gradient update to obtain the Instantaneous Control algorithm, see Algorithm~\ref{ICAlg}.
\RestyleAlgo{boxruled}
\begin{algorithm2e}[h!]
\caption{Instantaneous Control Algorithm}\label{ICAlg}
 \KwData{Initial data of $y$ or $p$; parameter values}
 \KwResult{Instantaneous control $\con$ or $\textbf{w}$; the corresponding states $y(\con)$ or $p(\textbf{w})$; optimal functional values}
 initialization\;
 $t=0; dt=T/K$\;
 \While{$t < T$}{
 \begin{itemize}
 \item [0.] solve state system \eqref{ODEsys} (or \eqref{PDEsys}, resp.)\;
 \item [1.] solve adjoint problem given in \eqref{rescaled_ad} (or \eqref{adSys}, resp.)\;
 \item [2.] compute gradient corresponding to \eqref{gradient}\;
 \item [3.] compute step size using the Armijo rule \eqref{ArmRule}
 \item [4.] update controls by steepest descent step \eqref{eq:update}\;
 \item [5.] project control onto the feasible set using \eqref{projection}\;
 \item [6.] initialize controls for the next time slice corresponding to \eqref{IC_update_control}\;
 \end{itemize}
  t = t+dt
 }
\end{algorithm2e}
\begin{rem}
 Note that we perform only one gradient step on each time slice. Numerical tests have shown that more gradient steps would only marginally improve the results, but drastically increase the computational effort. Further, we use the time slices of the instantaneous approach as time discretization for the simulation. There are no intermediate steps computed on one time slice.
\end{rem}

\section{Numerical Results}\label{sec:NumRes}
The numerical simulations for the mean-field equation are performed on the computational domain  $ \Omega = [-100,100]^2 \times [-5,5]^2\subset \mathbb{R}^{4}$, i.e., we set $D=2$. Further, we use the scaled CFL condition
\begin{equation}
 \frac{\tau |V|T}{Lh} \le 0.5
\end{equation}
with $L = 200, |V| = 5, dt = 0.02$ and thus $K = T/dt = 500.$

The grid parameter $h$ is varied throughout the simulations to investigate the grid convergence of the scheme. In fact, we use $25,50$ or $100$ grid points in each of the two directions leading to $h = 0.04, 0.02$ or $0.01$.  

Our particular choice of the interaction potentials are the Morse potentials as proposed in  \cite{Dorsogna, CarrilloMatrinPanferov}. For fixed positive parameters  $A_j,a_j,R_j,r_j$ we have 
\begin{equation}\label{InteractionPotentials}
 \Phi_j(x-y) = R_j \exp\left(-\frac{|x-y|}{r_j}\right) - A_j \exp\left( - \frac{|x-y|}{a_j}\right), \quad j=1,2.
\end{equation}
The parameters $A_j,R_j$ denote the attraction and repulsion strengths and $a_j, r_j$ the radii of interaction. The case $j=1$ refers to the interaction of the individuals, the interaction of individuals with external agents is denoted by $j=2$. Inspired by \cite{NoScaling}, we use the values
\begin{equation*}
 A_1 = 20,\quad R_1=50,\quad a_1=100,\quad r_1=2, \qquad\qquad A_2 = 5,\quad R_2=100,\quad a_2=1000,\quad r_2=50.
\end{equation*}

\begin{rem}
Note, that the coefficients of the interaction potentials are scaled with respect to the spatial size of the domain $\Omega$. As mentioned in the introduction we had the model problem of dogs herding sheep to the destination $\Edes$ in mind. The parameters for the sheep-sheep interaction have long-range attraction and repulsion on a very short range. The parameters modeling the sheep-dog interaction have a larger repulsive influence in order to reflect the guiding property correctly. 
\end{rem}

\begin{figure}[htbp]
	\begin{minipage}{0.49\textwidth} 
	\includegraphics[width=1.\textwidth]{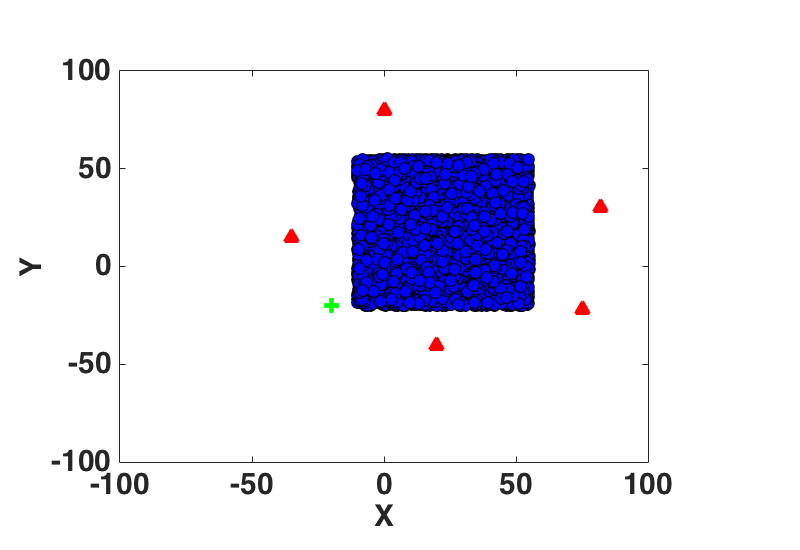}
	\end{minipage}
	\hfill
	\begin{minipage}{0.49\textwidth}
	\includegraphics[width=1.\textwidth]{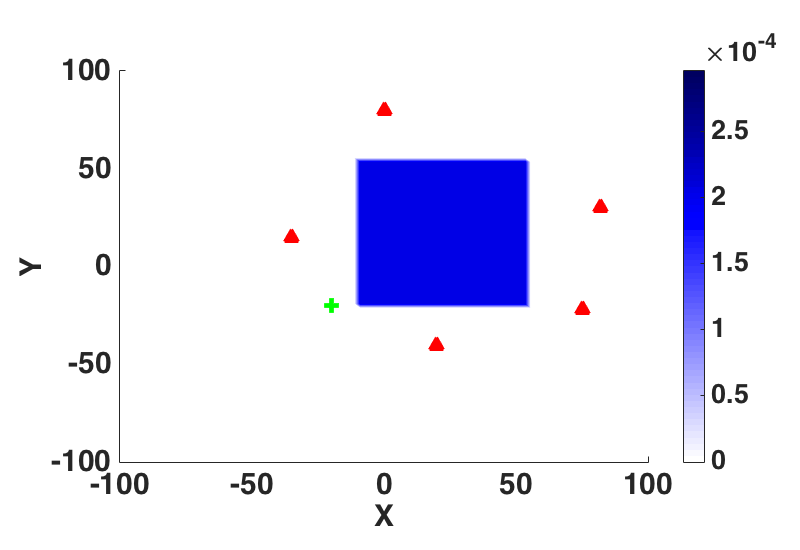}
	\end{minipage}
\caption{Initial configurations. Left: microscopic case with 8000 sheep represented by blue dots. Right: initial probability distribution $f$ for the mean-field case. In both cases the positions of the dogs $d$ are marked by red triangles and the destination $\Edes$ by a green dot.}
\label{fig:initial}
\end{figure}
Further parameters are fixed as follows: the destination $\Edes = (-20,-20)$ is depicted in green. The initial spatial support of the sheep is set to $\Omega_0 = [-10,55]\times[-20,55]$, i.e., $f_0(x,v)$ is the uniform distribution and the initial positions and velocities of the sheep are chosen by realizations of i.i.d.~random variables with $\text{law}(f_0)$. In every figure, the sheep and dogs are represented by blue markers and red triangles, respectively. The initial configurations of the microscopic and the mean-field case are shown in Figure~\ref{fig:initial}. 

The parameters for the following calculations are set to
\begin{equation*}
 T=10,\quad \bar{V}_n = 0.9\,\mathbb{V}_N(\x_0),\quad \bar{V}_\infty = 0.9\,\mathbb{V}_\infty(f_0) \quad \text{and} \quad \sigma_3 = 10^{-7}.
\end{equation*}
Thus, the desired variance is $10\%$ less than the initial variance given by the initial distribution of the crowd. 
The Armijo parameter is set to
$
 \omega_0 = 1\cdot 10^{3}.
$
\begin{rem}
Note, that in the given time interval $[0,T]$ the task of steering the crowd to the destination $\Edes$ is impossible to realize. In fact, our focus lies on the comparison of the behavior of the controls as $N$ increases. 
\end{rem}

\subsection{Numerical Results Using the IC Algorithm}
\begin{figure}
 	\begin{minipage}{0.49\textwidth}
	\includegraphics[width=1.\textwidth]{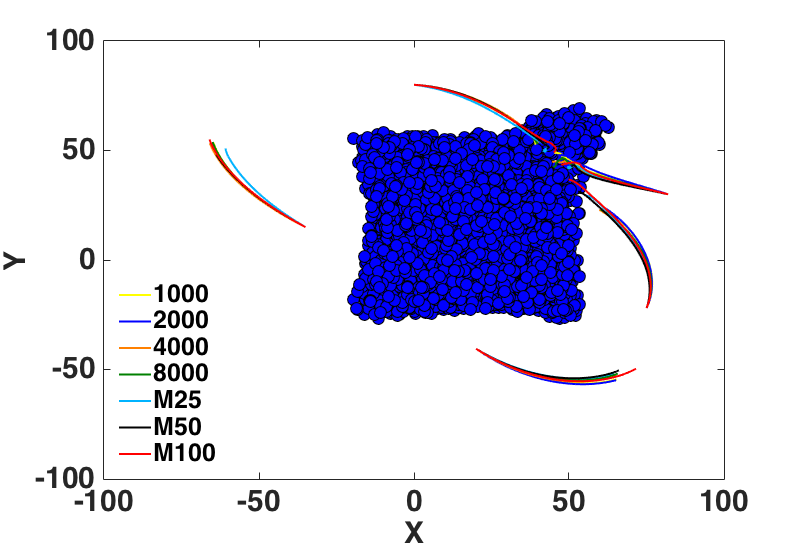}
	\end{minipage}
	\hfill
	\begin{minipage}{0.49\textwidth}
	\includegraphics[width=1.\textwidth]{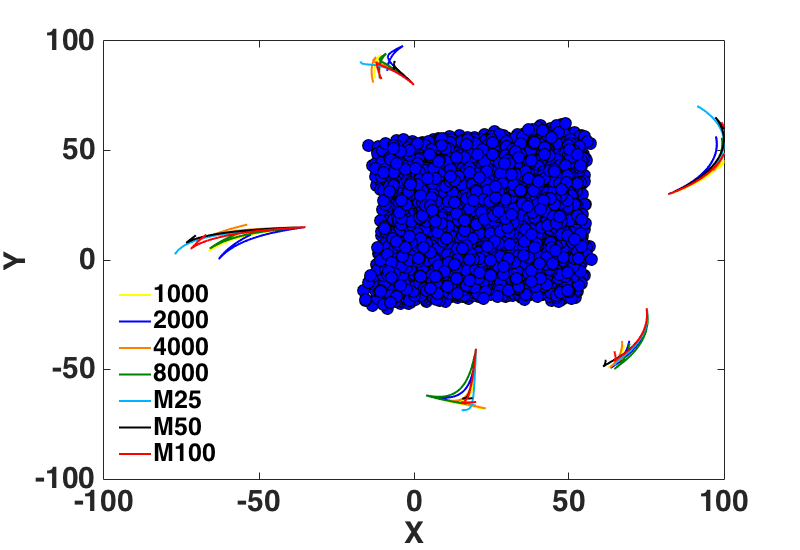}
	\end{minipage}
	\caption{Influence of the cost functional parameters. Left: Setting IC\ref{S1} - the dogs positioned on the left side and the bottom of the crowd clear the way to the destination. The others push the crowd from behind towards the destination. Right: Setting IC\ref{S2} - the dogs keep distance and orient towards to corners to reduce the variance of the crowd. These observations perfectly agree with the intention we had when modeling the cost functional.}
	 \label{fig:influence}
\end{figure}

In this section we discuss the results obtained by the IC algorithm. First, we illustrate the influence of the cost functional parameters. Therefore, we set up three test cases: \ref{S1} stresses the variance part $J_1$ of the cost functional, we expect the dogs to move towards the corners of the crowd in order to reduce the spread of the crowd. In $\ref{S2}$ the second part $J_2$ is emphasized. Thus, the dogs would be more inclined to push the crowd towards the destination. For test case \ref{S3} we choose the weights of the cost functional such that the focus lies on steering the crowd to the destination $\Edes$ while the variance term has minor influence but cannot be neglected. In particular, we use the following parameter values
\begin{align}\label{S1} \tag{\textbf{S1}}
 \sigma_1 &= 9\cdot 10^{-2}, &&\sigma_2 = 10^{-3}, \\
 \label{S2} \tag{\textbf{S2}}
 \sigma_1 &= 10^{-4}, &&\sigma_2 = 9\cdot 10^{-1}, \\
 \label{S3} \tag{\textbf{S3}}
 \sigma_1 &= 5\cdot 10^{-3}, &&\sigma_2 = 5\cdot 10^{-1}. 
\end{align}
These choices assure that $J_1$ of test case \ref{S1} has the same order of magnitude as $J_2$ of test case \ref{S2} and the other way around. The following notation is used in the Figures:
\begin{equation*}
 J_1 = \frac{\sigma_1}{4} \Big(\mathbb{V}_N(\x_t) - \bar{V}_N\Big)^2 \quad\text{and}\quad J_2 = \frac{\sigma_2}{2} \norm{\mathbb{E}_N(\x_t)-\Edes}_2^2.
\end{equation*}
 The spatial and the velocity grid are discretized with the same number of grid points ($25, 50$ or $100$) in each of the two directions. 
The graphs corresponding to mean-field solutions are labeled M$\#$ where $\#$ denotes the number of grid points. The graphs corresponding to microscopic simulations are denoted by their respective number of particles. In Figure~\ref{fig:influence} the influence of the cost functional parameters is illustrated. 

\begin{figure}
	\begin{minipage}{0.49\textwidth}
	\includegraphics[width=1.\textwidth]{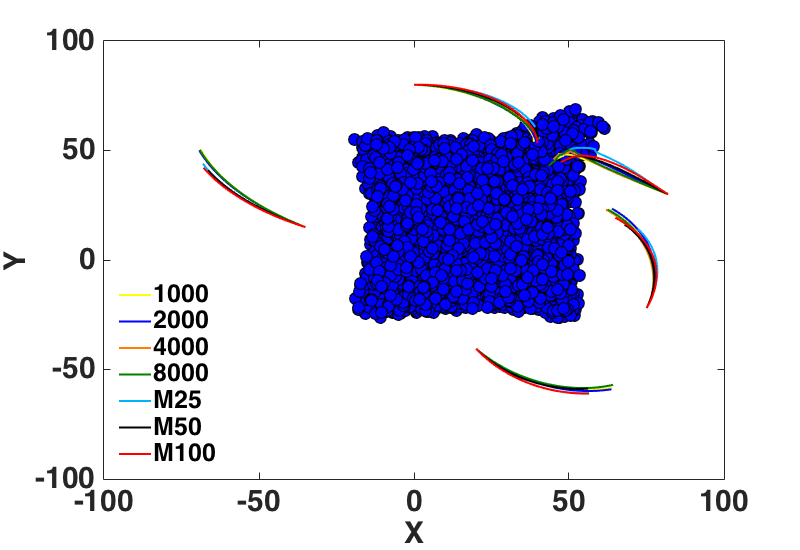}
	\end{minipage}
	\hfill
	\begin{minipage}{0.49\textwidth}
	\includegraphics[width=1.\textwidth]{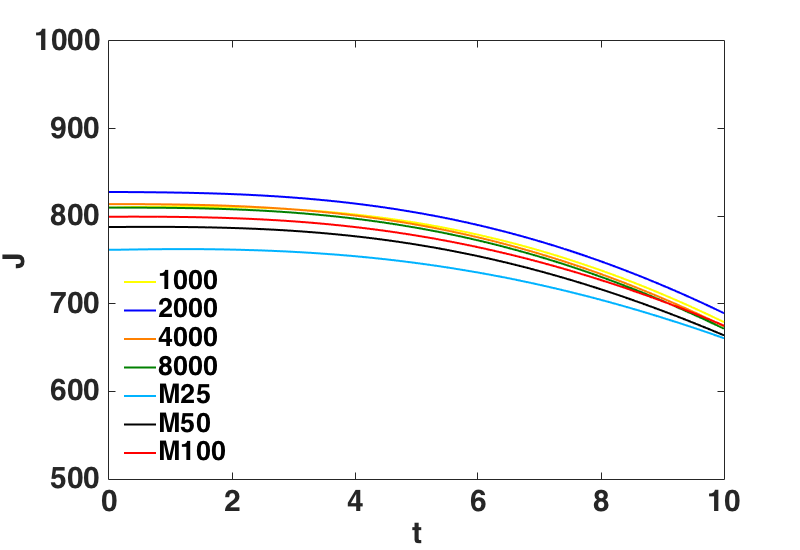}
	\end{minipage}
	\caption{Setting IC\ref{S3} - Trajectories of the dogs and evolution of the cost functional. Left:  the dogs positioned on the left side and at the bottom of the crowd clear the way to the destination. The others push from behind. Note that the dogs on the left to not enter the crowd as deep as in setting IC\ref{S1}. Right: Evolution of the cost functional values.}
	\label{fig:pathAndCost}
\end{figure}

Considering the application we expect that \ref{S3} is the most realistic setting, since the focus lies on $J_2$ which measures the distance to the destination, while the influence of $J_1$ is not negligible. In Figure~\ref{fig:pathAndCost}(left) the trajectories of the dogs and the crowd at $T=10$ are shown. Note, that the dogs are not splitting the crowd as much as in test case \ref{S1}. 

\begin{figure}
	\begin{minipage}{0.49\textwidth}
	\includegraphics[width=1.\textwidth]{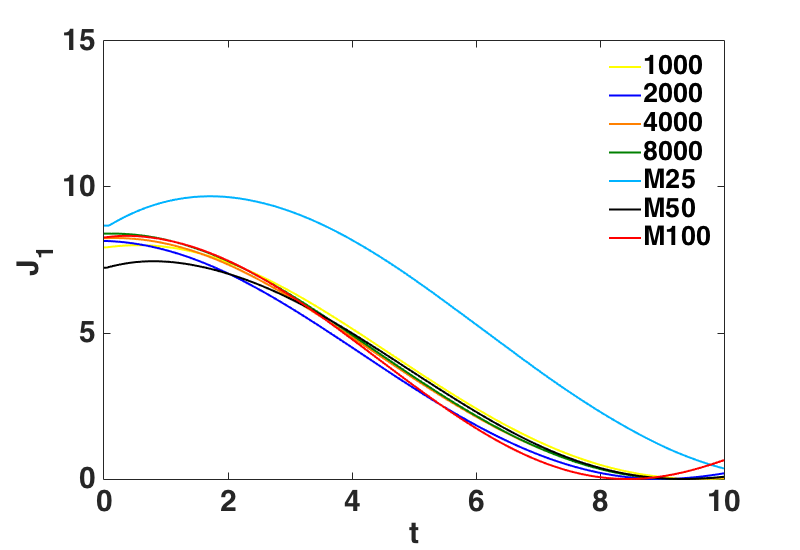}
	\end{minipage}
	\hfill
	\begin{minipage}{0.49\textwidth}
	\includegraphics[width=1.\textwidth]{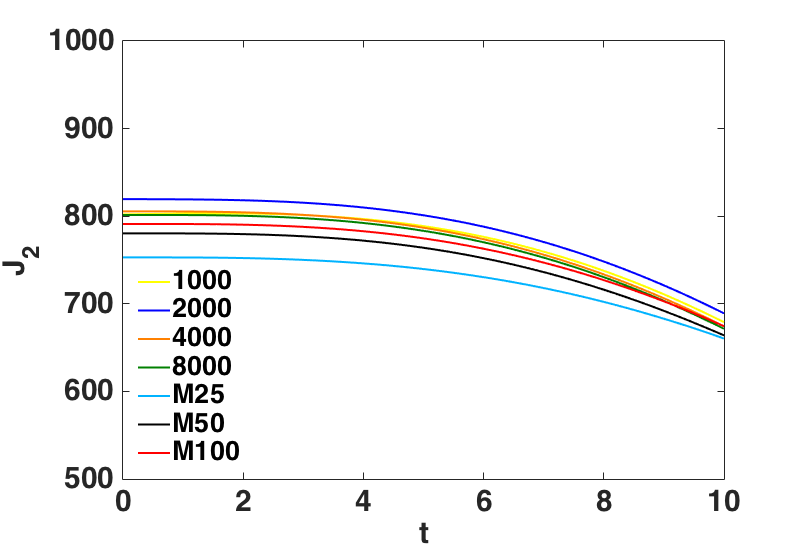}
	\end{minipage}
	 \caption{Setting IC\ref{S3} - Evolution of the cost functional parts. Left: values of $J_1$. The simulation $M25$ overestimates the variance significantly. All other graphs are in good agreement. Right: values of $J_2$. All graphs show similar behavior. We see convergence as $N\rightarrow \infty$ and as $h \rightarrow 0$.}
	 \label{fig:costDetail}
\end{figure}

Comparing the values of $J$, $J_1$ and $J_2$ in Figure~\ref{fig:pathAndCost}(right) and Figure~\ref{fig:costDetail} we conclude that the part measuring the crowd's distance to the destination is dominating the others. 
\begin{rem}
Note, that the graphs of the cost functional values are given with respect to time, i.e., it may happen that the values increase. Since the initial positions of the dogs is chosen arbitrarily, they move to appropriate positions first, causing a slight increase of the cost functional.  Afterwards the cost functional is decreasing as expected. 
\end{rem}

The simulation on the coarse mean-field discretization $M25$ overestimates the variance term significantly as can be seen in Figure~\ref{fig:costDetail}(left). The evolution of $J_1$ is in good resemblance for all other discretizations. Similarly, there is a good agreement  in the graphs showing the evolution of $J_2$. Furthermore, we see the convergence of the graphs to the solution for $M100$ as $N$ increases and the convergence of the mean-field simulations as the grid is refined.

In Table~\ref{Table_IC} we illustrate the convergence for $N \rightarrow \infty$ by the  relative errors of
\begin{subequations} \label{num_norms}
\begin{align}
\norm{u - u_\text{ref}} &= \int_0^T \norm{u(t) - u_\text{ref}(t)}_{ \mathbb{R}^4} \dd t, \\
|J- J_\text{ref}|&=  \int_0^T |J(t) -J_\text{ref}(t) |  \dd t, \\
|\rho^N - \rho_\text{ref}| &= \int \norm{\rho(t,x,y) - \rho_\text{ref}(t,x,y)} \dd x \dd y\dd t.
\end{align}
\end{subequations}
The reference values are the results of the simulation $M100$. To compute the norms the empirical density $\rho^N$ of the microscopic simulations is approximated by a histogram which is based on the grid of the corresponding mean-field simulation. 

In the first two columns we see the convergence of the mean-field scheme, as expected the values are decreasing as the grid is refined. The four columns on the right illustrate the behavior as $N$ increases. All norm values are decreasing for increasing number of particles from $N=1000$ to $N=8000$. 

\begin{rem}\label{integral_quantities}
Since mean-field quantities are averaged it is very common to compare integral values like $J$ and $\rho^N$ in this setting. Note that the velocities and the trajectories of the agents are no such quantities. This explains why the convergence is more transparent for the integral quantities than for the velocities. 
\end{rem}

\begin{figure}
 	\begin{minipage}{0.49\textwidth}
	\includegraphics[width=1.\textwidth]{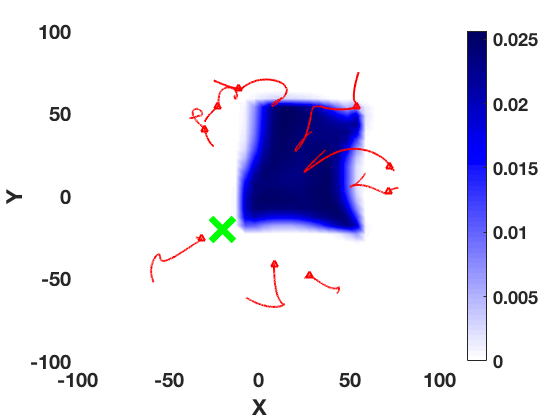}
	\end{minipage} 
	\hfill
	\begin{minipage}{0.49\textwidth}
	\includegraphics[width=1.\textwidth]{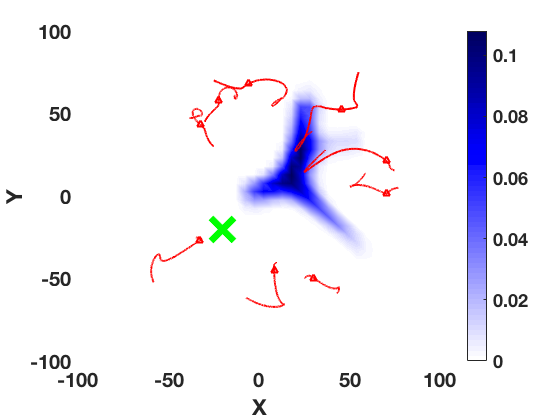}
	\end{minipage}
	\begin{center} $t=3$ \hspace{7cm} $t=9$ \end{center}
 	\begin{minipage}{0.49\textwidth}
	\includegraphics[width=1.\textwidth]{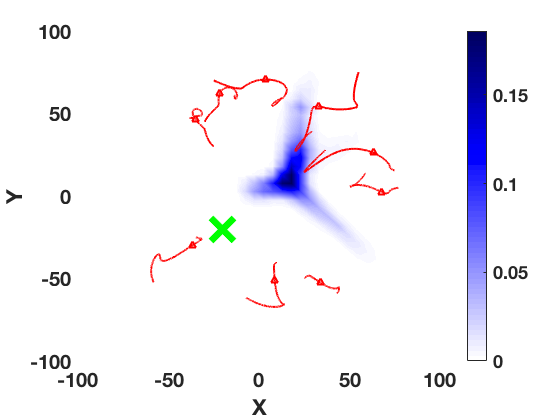}
	\end{minipage}
	\hfill
	\begin{minipage}{0.49\textwidth}
	\includegraphics[width=1.\textwidth]{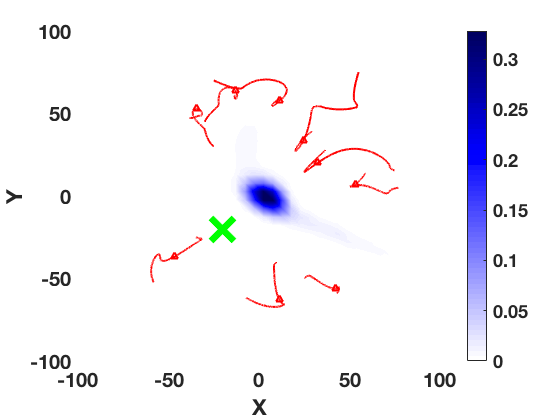}

	\end{minipage}
	\begin{center} $t=12$ \hspace{7cm} $t=20$\end{center}
 	\begin{minipage}{0.49\textwidth}
	\includegraphics[width=1.\textwidth]{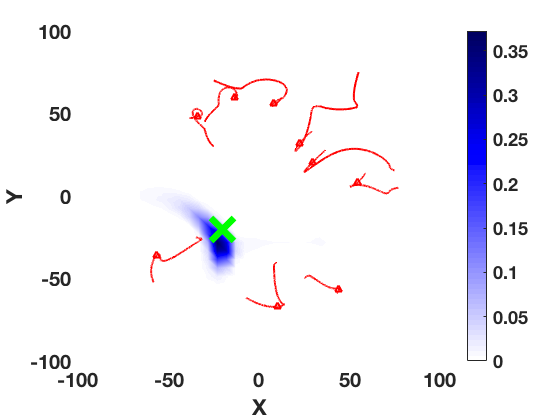}
	\end{minipage}
	\hfill
	\begin{minipage}{0.49\textwidth}
	\includegraphics[width=1.\textwidth]{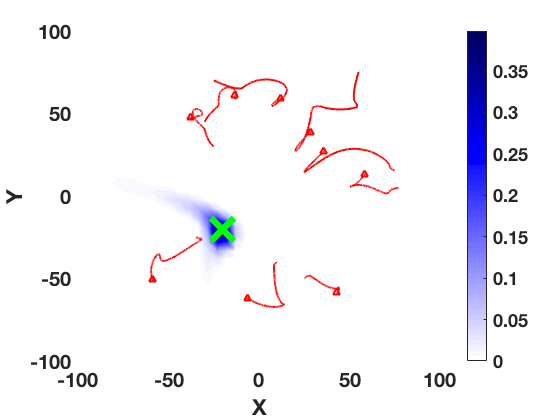}
	\end{minipage}
	\begin{center} $t=40$ \hspace{7cm} $T=55$\end{center}
	\caption{ Mean-field optimization for a long time horizon $T=55$ with 9 dogs. The destination is marked in green. The trajectories of the dogs are depicted in red. The current position of each dog is marked by a red triangle. The dogs are guiding the crowd. First the variance of the crowd is reduced, then the center of mass is pushed towards the destination.}
	\label{fig:ParticleMFVergleich}
\end{figure}

\begin{table}[!ht]
\begin{tabular}{| c | c | c || c | c | c | c | }
 \hline
  IC\ref{S1}       &  M25    & M50    & 1000     & 2000  & 4000  & 8000         \\
 \hline
 $|J - J_\text{ref}| /|J_\text{ref}|$  \ &    1.2 & 0.2  & 0.9   & 1.1  & 0.8 & 0.75   \\
 \hline
 $\norm{ u- u_\text{ref}} / \norm{u_\text{ref}}$ &  4.21   & 1.51 & 1.33  & 1.33   &  1.32  & 1.32     \\
  \hline
 $|\rho^N - \rho_{\text{ref}}| / |\rho_\text{ref}|$ &  - & - & 0.1   & 0.1   & 0.1  &   0.1       \\
  \hline
  IC\ref{S2}     &      &     &      &   &  &        \\
 \hline
 $|J - J_\text{ref}| /|J_\text{ref}|$ &  0.84 & 0.1 & 0.1  &  0.04   & 0.05 & 0.05    \\
 \hline
 $\norm{ u- u_\text{ref}} / \norm{u_\text{ref}}$&  3.43   & 1.04  & 1.29  & 1.35   &  1.3 & 1.32    \\
  \hline
 $|\rho^N - \rho_{\text{ref}}| / |\rho_\text{ref}|$&  - & - & 0.01   & 0.01   & 0.1  &   0.01      \\
  \hline
   IC\ref{S3}     &      &     &      &   &  &          \\
 \hline
  $|J - J_\text{ref}| /|J_\text{ref}|$ &  0.68 & 0.12 & 0.14  &  0.29   & 0.13  & 0.09     \\
 \hline
 $\norm{ u- u_\text{ref}} / \norm{u_\text{ref}}$ & 2.8   & 1.01  & 1.21  & 1.21  &  1.22 & 1.21     \\
  \hline
$|\rho^N - \rho_{\text{ref}}| / |\rho_\text{ref}|$ &  - & - & 0.01   & 0.01   & 0.01  &   0.01      \\
  \hline
\end{tabular}
 \caption{Illustration of convergence as $N\rightarrow \infty$ using the norms in \eqref{num_norms} with the respective values of $M100$ as reference. All values are given in \%.}
 \label{Table_IC}
\end{table}
In Figure~\ref{fig:ParticleMFVergleich} the evolution of the crowd and the dogs is illustrated for the mean-field case $M50$ on the time interval $I = [0,55]$ with $\bar V_\infty = 0.1 \mathbb{V}_\infty(f_0)$. The red lines show the trajectories of the dogs as above. Additionally, to the information in Figure~\ref{fig:pathAndCost}, the red markers indicate the current positions of the dogs at different times. We see that  the variance of the crowd is reduced first and then the center of mass is transported to the destination.
\newpage
\begin{table}
	\begin{tabular}{| c | c | c || c | c | c | c | }
		\hline
		computation time     &  M25    & M100    & 1000     & 2000  & 4000  & 8000         \\
		\hline
		$t /t_\text{ref}$   &    0.15 & 13  & 0.16   & 0.64  & 2.74 & 10.75        \\
		\hline \hline
		memory     &  M25    & M100    & 1000     & 2000  & 4000  & 8000    \\
		\hline
		$m/m_\text{ref}$   &    0.08 &   15.84   & 0.15  & 0.56 & 2.20 & 9.26   \\
		\hline
	\end{tabular}
	\caption{Illustration of relative computational times and memory consumption for 10 time steps of the optimization algorithm. The reference values are the values of M50, respectively.}
	\label{tab:computationTimes}
\end{table}
In Table~\ref{tab:computationTimes} we compare the relative computation times and memory needed to compute 10 time steps of the optimization procedure. Replacing the 8000 particle simulation with the mean-field simulation with 50 grid points in each direction reduces the computational time by factor 10. Similar savings can be expected regarding the memory consumption.

\section{Conclusions}
Altogether, the numerical tests with different cost functional parameters show that the cost functional yields in the optimization the expected behaviour. Further, the results underline the convergence of the states and the convergence of the controls as $N$ tends to infinity. The convergence will be investigated from an analytical point of view in future work. Finally, comparing the first order necessary conditions, one realizes that there appear derivatives of $S$ on the microscopic level but not in the corresponding the mean-field equations. This relation will be discussed in some upcoming  work as well.

\section*{Acknowledgements} 
MB acknowledges support by ERC via Grant EU FP 7 - ERC Consolidator Grant 615216 LifeInverse.

\bibliographystyle{plain}
\bibliography{biblio}
\end{document}